\documentclass[12pt, reqno]{amsart}
\usepackage{ifthen,amsfonts,amsmath,amssymb,epic,eepic}
\usepackage{epsfig,color}

\makeatother

\theoremstyle{definition}

\newtheorem{Thm}{Theorem}
\newtheorem{Lem}{Lemma}

\newtheorem{Pro}{Proposition}

\numberwithin{equation}{section}

\newcommand{\nn}{{\bf{n}}}

\def\RR{{\bold R}}

\newcommand{\dv}{{\text {div}}}
\newcommand{\e}{{\text {e}}}
\newcommand{\Area}{{\text {Area}}}

\newcommand{\eqr}[1]{(\ref{#1})}

\begin{document}

\title[Sharp estimates for mean curvature flow of graphs]
{Sharp estimates for mean curvature flow of graphs}

\author{Tobias H. Colding}%
\address{Courant Institute of Mathematical Sciences\\
251 Mercer Street\\ New York, NY 10012}
\author{William P. Minicozzi II}%
\address{Department of Mathematics\\
Johns Hopkins University\\
3400 N. Charles St.\\
Baltimore, MD 21218}
\thanks{The authors were partially supported by NSF Grants DMS
0104453 and DMS 0104187}


\email{colding@cims.nyu.edu, minicozz@jhu.edu}

\maketitle

\section{Introduction}

A one-parameter family of smooth hypersurfaces $\{ M_t \} \subset
\RR^{n+1}$ {\it flows by mean curvature} if
\begin{equation}
    z_t = {\bf{H}} (z) = \Delta_{M_t} z \, ,
\end{equation}
where $z$ are coordinates on $\RR^{n+1}$ and ${\bf{H}} = - H \nn$
is the mean curvature vector.

In this note, we prove  sharp gradient and area estimates for
graphs flowing by mean curvature.  Thus, each $M_t$ is assumed to
be the graph of a function $u(\cdot , t)$.  So, if $z = (x,y)$
with $x \in \RR^n$, then $M_t$ is given by $y = u (x,t)$.
   Below, $du$ is the $\RR^n$
gradient of a function $u$, $\|u\|_{\infty}$ is the sup norm, and
$B_{s}$ is the ball in $\RR^n$ with radius $s$ centered at the
origin.

 Our gradient estimate is the
following (see Section \ref{s:area} for the sharp area estimate):

\begin{Thm}     \label{t:gbp}
There exists $C= C(n)$ so if the graph of $u: B_{\sqrt{2n+1} r}
 \times [0,r^2] \to \RR$ flows by mean
curvature, then
\begin{equation}    \label{e:gbp}
    \log |du|(0,r^2/[4n]) \leq   C \, ( 1 + r^{-1} \,
\|u(\cdot , 0) \|_{\infty} )^2 \, .
\end{equation}
\end{Thm}

The quadratic dependence on  $\|u(\cdot , 0) \|_{\infty}$ in
\eqr{e:gbp} should be compared with the linear dependence which
holds when the graph of $u$ is minimal (i.e., $u_t = 0$).   In the
minimal case, Bombieri, De Giorgi, and Miranda proved in
 \cite{BDM} that
\begin{equation}
    \log |du|(0) \leq   C \, ( 1 + r^{-1} \,
    \|u \|_{\infty} )
\end{equation}
(the case of surfaces was done by Finn in \cite{F1}).  By an
earlier example  of Finn, this exponential dependence cannot be
improved even in the minimal case (see \cite{F2} and cf.
\cite{GiTr}).

In \cite{K}, Korevaar gave a maximum principle proof of a weaker
form of \cite{BDM}; this weaker form had $\|u\|^2_{\infty}$ in
place of $\|u\|_{\infty}$. Ecker and Huisken adapted Korevaar's
argument to mean curvature flow in theorem 2.3 of \cite{EH2} to
get
\begin{equation}    \label{e:eh2}
    \log |du|(0,r^2/[4n]) \leq  1/2 \, \log
    \left( 1 + \|du (\cdot,0)\|^2_{\infty} \right)
      + C \, ( 1 + r^{-1} \, \|u(\cdot , 0) \|_{\infty} )^2 \, .
\end{equation}
Note that, unlike \eqr{e:gbp}, the gradient bound \eqr{e:eh2}
depends also on the initial bound for the gradient.

\vskip2mm Using the so-called grim reaper, one can see that the
quadratic dependence on $\| u (\cdot , 0)  \|_{\infty}$ in Theorem
\ref{t:gbp} is sharp (see Proposition \ref{p:grr} below). The {\it
grim reaper} is the translating solution to the mean curvature
flow given by that for each $t$ it is a graph of the function
\begin{equation}    \label{e:gr}
    u(x,t) = t - \log \sin x \, ,
\end{equation}
where $x  \in (0 , \pi )$ and $t \in [0,\infty)$.  Note that
$-u(x,t)$ is a downward translating solution.   More generally, a
parabolic rescaling by $\lambda>0$ gives that the graph of
\begin{equation}
u^{\lambda}(x,t)=\frac{1}{\lambda}u(\lambda  x ,\lambda^2 t) =
    \lambda \, t - \frac{\log \sin (\lambda x )}{\lambda}  \,  ,
\end{equation}
where $x  \in (0 , \pi / \lambda )$, is a translating solution
 flowing with speed $\lambda$; see figure \ref{f:f1}.  
Since $\lim_{x\to 0} \frac{\sin x}{x} = 1$,
 an easy calculation shows that  for
  $\lambda > 0$ sufficiently large
\begin{equation}    \label{e:ulam}
    u^{\lambda} ( \e^{-\lambda^2} , 1) = \lambda
    - \frac{\log \sin (\lambda \e^{-\lambda^2} )}{\lambda} \leq  2
    \lambda \, .
\end{equation}

\begin{figure}[htbp]
    \setlength{\captionindent}{4pt}
    \begin{minipage}[t]{0.5\textwidth}
    \centering\input{grad1.pstex_t}
    \caption{Two scaled grim reapers.}\label{f:f1}
    \end{minipage}\begin{minipage}[t]{0.5\textwidth}
    \centering\input{grad2.pstex_t}
    \caption{Use $u^+ = u^{\lambda}(x + \frac{\pi}{\lambda} , t) - 3 \lambda$ 
     and $u^-= -u^{\lambda}(x,t) + 3 \lambda$ as barriers and let
    $w$ be a graph between $u^+$ and $u^-$.}\label{f:f2}
    \end{minipage}
\end{figure}

\begin{Pro}     \label{p:grr}
Given $\lambda > 1$ sufficiently large, there is a solution
$w(x,t)$ on $\RR \times [0,\infty)$ of the mean curvature flow
with
\begin{align}   \label{e:grr1}
    3 \lambda &< \|w (\cdot , 0 ) \|_{\infty} \leq 4 \lambda \, , \\
    \lambda \e^{\lambda^2} & \leq \max_{ |x| \leq \e^{-\lambda^2} } |d w (x,1) |   \, .
    \label{e:grr2}
\end{align}
\end{Pro}

\begin{figure}[htbp]
    \setlength{\captionindent}{4pt}
    \begin{minipage}[t]{0.5\textwidth}
    \centering\input{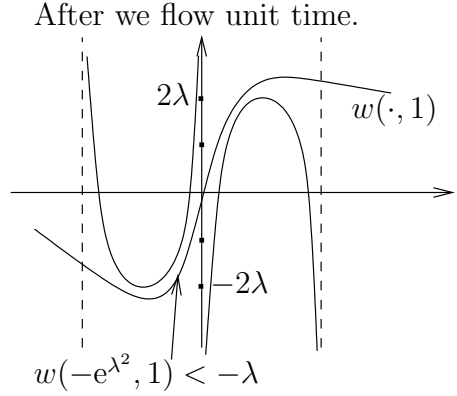}
    \caption{By the maximum principle, $w(\cdot , 1)$ will be between the
    barriers. This bounds $\max |\nabla w|$ from below.}\label{f:f3}
    \end{minipage}
\end{figure}

\begin{proof}
Define solutions $u^{+} (x,t)= u^{\lambda}(x + \pi/\lambda ,t) - 3
\lambda$ for $-\pi/\lambda < x < 0$
 and $u^{-}(x,t) = -u^{\lambda}(x,t) + 3\lambda$
for $0<x<\pi/\lambda$ of the mean curvature flow to be used as
barriers. Since  $u^{+} \geq -3\lambda $ and $u^{-} \leq
3\lambda$, it is easy to choose (see figure \ref{f:f2})
a smooth compactly supported
function $w(\cdot , 0): \RR \to \RR$ satisfying \eqr{e:grr1} and
so
\begin{align}
    w(x,0) &< u^{+} (x,0) {\text{ for }} -\pi/\lambda < x < 0
    \, , \label{e:up} \\
    u^{-}(x,0) &< w (x,0) {\text{ for }} 0 < x < \pi/\lambda
    \, . \label{e:down}
\end{align}
(We can choose $w(\cdot , 0)$ so that $w(x,0) = 0$ for $|x|
>  \pi / \lambda$.)  The existence results of \cite{EH1}
or \cite{EH2} (see, e.g., theorem 1.7 in \cite{E}) extend $w(x,0)$
to a solution $w(x,t)$ of the mean curvature flow defined for $x
\in \RR$ and $t \in [0, \infty)$; see figure \ref{f:f3}. 
Moreover, the maximum principle
extends \eqr{e:up} and \eqr{e:down} to all $t \geq 0$.  In
particular, using this at $x = \pm \e^{-\lambda^2}$, $t=1$ and
substituting \eqr{e:ulam}, we get that
\begin{align}        \label{e:punchline1}
    w(-\e^{-\lambda^2} , 1) &< u^{+} (-\e^{-\lambda^2} , 1) \leq - \lambda \, , \\
     \lambda
    &\leq     u^{-} (\e^{-\lambda^2} , 1)          < w(\e^{-\lambda^2} , 1) \,
    .  \label{e:punchline2}
\end{align}
Finally, combining \eqr{e:punchline1}, \eqr{e:punchline2}, and the
mean value theorem gives \eqr{e:grr2}; see figure \ref{f:f3}.
\end{proof}

\vskip2mm
 Throughout,
 $\nabla$, $\Delta$, and $\nn$ are
the induced covariant derivative, laplacian, and unit normal on
the submanifold $M_t$  of $\RR^{n+1}$.  The  graph of a function
$u$ flows by mean curvature if
\begin{equation}    \label{e:evq}
    u_t =  (1 + |du|^2)^{1/2} \,
    \dv \left( \frac{du}{(1 + |du|^2)^{1/2} } \right)
    \, ,
\end{equation}
where $\dv$ is divergence in $\RR^n$.

\section{Interior gradient estimate for graphs}

Theorem \ref{t:gbp}   will follow immediately from the next
proposition and a standard maximum principle bounding $u$ at
future times in terms of the initial bound, see Lemma
\ref{l:sphere}.

\begin{Pro}     \label{p:gbp}
If  the graph of $u: B_r \times [0,r^2] \to \RR$ flows by mean
curvature, then
\begin{equation}    \label{e:gbp2}
    \log \left( 1 + |du|^2 (0,r^2/[4n]) \right)
        \leq 2  \log 10
+ 16 \, n \, (1 + 2 \, r^{-1} \, \|u  \|_{\infty}  )^2  \, .
\end{equation}
\end{Pro}

The strategy of the proof of Proposition \ref{p:gbp} is as
follows.  By the maximum principle (Lemma \ref{l:mp2}) at the
maximum of $\phi\,v$ (where $\phi$ is a cutoff function and $v=
\left( 1 + |du|^2 \right)^{1/2}$) the heat operator of the cutoff
function is nonnegative. By choosing an appropriate cutoff
function in terms of $u$ (Lemma \ref{l:cop}), we can bound the
heat operator of the cutoff from above in terms of the gradient of
$u$. Playing off this lower and upper bound at the max against
each other  gives the proposition.

\begin{Lem}     \label{l:mp2}
Suppose $( \partial_t - \Delta  ) v \leq - 2 |\nabla v|^2 / v$ for
a  function $v \geq 0$ on $\{ M_t \}_{t \in [0,1]}$.  If the
function $\phi$ is $\leq 0$ on  $M_0 \cup_t (\partial M_t)$ but
$\max \phi v > 0$, then   $(\partial_t - \Delta ) \phi \geq 0$ at
the maximum of $\phi v$.
\end{Lem}

\begin{proof}
At the maximum of $\phi v$, we get
\begin{align}
   \nabla (\phi v) &= v  \nabla \phi + \phi \nabla v = 0 \, , \label{e:mp1} \\
   (\partial_t - \Delta) (\phi v) &=  v \partial_t \phi +
     \phi \partial_t v -
    v \Delta \phi - 2 \langle \nabla \phi ,
   \nabla v \rangle - \phi \Delta v \geq 0 \, . \label{e:mp2}
\end{align}
Substituting \eqr{e:mp1} into \eqr{e:mp2} and using $( \partial_t
- \Delta  ) v \leq - 2 |\nabla v|^2 / v$ gives $(\partial_t -
\Delta ) \phi \geq 0$.
\end{proof}

We will apply Lemma \ref{l:mp2} to the volume element $v = (1 +
|du|^2)^{1/2}$. We will first need some elementary formulas. If
the graph of $u$ flows by mean curvature and $a \in \RR$, then
\begin{equation}
       ( \partial_t - \Delta  ) \,
          \e^{ay^2/t}  = - \frac{\e^{ay^2/t}}{t^2} \,
       \left[ a  y^2  + 4 a^2 y^2 |\nabla y|^2    + 2a t
       |\nabla y|^2 \right]  \, , \label{e:y2}
\end{equation}
and (see, e.g., lemma 1.1 in \cite{EH2} or 2.11 in \cite{E})
\begin{align}
    ( \partial_t - \Delta  )  v = -|A|^2 \, v -
2 \frac{|\nabla v|^2}{v} &\leq - 2 \frac{|\nabla v|^2}{v} \,
    , \label{e:jac} \\
    ( \partial_t - \Delta  )  (1 -|x|^2-2nt) &\leq 0 \, . \label{e:dx2}
\end{align}

The next lemma introduces the cutoff $\phi$ which will be used in
Lemma \ref{l:mp2}.

\begin{Lem}         \label{l:cop}
Set $\phi =   \eta \, \e^{ay^2/t}$ for $a \in \RR$ and $\eta =
 (1 -|x|^2-2nt)$. If the graph of $u:B_1 \times [0,1] \to \RR$
flows by mean curvature,
 then
\begin{equation}     \label{e:co1p}
    (\partial_t - \Delta ) \phi \leq -
\frac{ \e^{ay^2/t} }{t^2} \, \left[ a y^2 \eta + (4 a^2 y^2 + 2at)
\eta \frac{|du|^2}{1 + |du|^2}
    - 8 |ay| \, t \frac{|du|}{1 + |du|^2}  \right]  \,
 \, .
\end{equation}
\end{Lem}

\begin{proof}
Using \eqr{e:y2} and \eqr{e:dx2} gives
\begin{align}        \label{e:co2p}
    (\partial_t - \Delta) \phi &=  \eta  (\partial_t - \Delta)
    \,   \e^{ay^2/t}  +  \e^{ay^2/t}
    (\partial_t - \Delta)  \eta - 4 \e^{ay^2/t} \frac{a y}{t}
\langle \nabla \eta ,
    \nabla y \rangle \notag \\
& \leq - \frac{\e^{ay^2/t} }{t^2} \left[ (a y^2 + 4 a^2 y^2
|\nabla y|^2 + 2at |\nabla y|^2) \eta  -  4 ay \, t \langle \nabla
y , \nabla |x|^2 \rangle \right] \, .
\end{align}
The lemma follows since $|x| \leq 1$ and the $y$ component of the
normal is $(1 + |du|^2)^{-1/2}$.
\end{proof}

\begin{proof}
(of Proposition \ref{p:gbp}.) By scaling, it suffices to prove the
proposition when $r=1$.  Set $\eta = (1-|x|^2 - 2nt)$ and $\phi =
\eta \e^{a y^2/t}$ for $a \leq -2$ to be chosen. After replacing
$u$ by $u + \|u\|_{\infty} + 1$ (i.e., translating), we can assume
that $u \geq 1$; in particular, $\phi$ vanishes when $t=0$.   If
the maximum of $\phi v$ for $t\in [0,1]$ is at $(x_0 , y_0 ,t_0)
\in B_1 \times \RR \times (0,1]$, then \eqr{e:jac} together with
Lemmas \ref{l:mp2} and \ref{l:cop} give
\begin{equation}        \label{e:q}
        a y_0^2 \eta  + (4 a^2 y_0^2 +
2at_0) \eta \frac{|du|^2(x_0,t_0) }{1 + |du|^2 (x_0,t_0)}
    - 8 |ay_0| \, t_0 \frac{|du|(x_0,t_0)}{1 + |du|^2(x_0,t_0) } \leq 0 \, .
\end{equation}
There are now two cases.  Namely, either
 \begin{equation} \label{e:case1}
    |ay_0| \eta
|du| (x_0 , t_0) < 8 \, ,
\end{equation}
or $|ay_0| \eta |du| (x_0 , t_0) \geq  8$; in the second case,
\eqr{e:q} (and $4 a^2 y_0^2 + 2at_0 > 2 a^2 y_0^2$) yields
\begin{equation}        \label{e:case2}
        a y_0^2 \eta  +  a^2 y_0^2  \eta \frac{|du|^2(x_0,t_0) }{1 + |du|^2 (x_0,t_0)}
     \leq 0 \, .
\end{equation}
Since $\eta \leq 1$, we get in either case that
 \begin{equation} \label{e:cases}
     \eta |du| (x_0 , t_0) \leq 4 \, ,
\end{equation}
 Since $\max_{[0,1]} \, (\phi v) = \phi v(x_0, y_0 , t_0)$ and $a < 0$, we get
\begin{equation}    \label{e:gb3p}
     \phi v = \eta \, \e^{a \, y^2/t} \, \left( 1 +  |du|^2 \right)^{1/2}
     \leq 5
         \, .
\end{equation}
\end{proof}

A standard barrier argument using shrinking spheres  bounds   the
future height by the initial height:

\begin{Lem}     \label{l:sphere}
If $\rho \geq (2n+1)^{1/2}$ and the graph of $u: B_{ \rho r }
\times [0,r^2] \to \RR$ flows by mean curvature, then
\begin{align}    \label{e:sphere}
    \max_{ B_r \times [0,r^2]} |u(x,t)|
&\leq  r \, \left[  \rho  - \left( \rho^2 - (2n+1)
\right)^{1/2}  \right] + \max_{B_{ \rho r } } |u (x,0)| \notag \\
 &\leq     \frac{(2n+1) \, r}{ \rho} + \max_{B_{ \rho r } } |u (x,0)| \, .
\end{align}
\end{Lem}

\begin{proof}
By scaling, it suffices to prove the lemma when $r=1$.
 Recall that
the one-parameter family $M_t$ of concentric spheres in
$\RR^{n+1}$ of radius $( \rho^2 - 2nt )^{1/2}$ centered at $x=0$,
$y= \rho + \max_{B_{ \rho } } u (x,0) + \epsilon$ is a solution to
mean curvature flow. For $\epsilon > 0$, $M_0$ does not intersect
the graph of $u(\cdot , 0)$. Applying the maximum principle and
letting $\epsilon \to 0$, we get that
\begin{equation}
    \max_{ B_1 \times [0,1]} u(x,t) \leq
\max_{B_{ \rho } } u (x,0) + \rho - \left( \rho^2 - (2n+1)
\right)^{1/2} \, .
\end{equation}
This, and a similar argument for  the minimum of $u$, gives
 \eqr{e:sphere}.
\end{proof}

\section{Area estimates for graphs flowing by mean curvature}
\label{s:area}

In this section, we prove an area bound for graphs flowing by mean
curvature which depends quadratically on the $L^{\infty}$ norm of
the initial height (integrating our gradient estimate gives an
exponential bound).  We also give an example showing that this is
sharp.

\begin{Thm}  \label{t:areabound}
There exists $C= C(n)$ so if the graph of $u: B_{\sqrt{2n+1} r}
 \times [0,r^2] \to \RR$ flows by mean
curvature, then
\begin{equation}    \label{e:abp}
    \Area ( u (B_{r/2} , r^2) )=\int_{ B_{r/2} }(1+|du|^2)^{1/2} (x,r^2) \, dx
\leq   C \,r^{n}\, \left( 1 + r^{-1} \, \|u(\cdot,0)\|_{\infty}
\right)^2 \, .
\end{equation}
\end{Thm}

Before proving  Theorem \ref{t:areabound}, we first argue as in
Proposition \ref{p:grr} to see that  the quadratic dependence on
$\|u(\cdot,0)\|_{\infty}$ is sharp:

\begin{Pro}     \label{p:grr2}
Given an integer $k > 1$, there is a solution $w(x,t)$ on $\RR
\times [0,\infty)$ of the mean curvature flow with
\begin{align}   \label{e:grr4}
    2 k &< \|w (\cdot , 0 ) \|_{\infty} \leq 3 k \, , \\
    4 k^2 - 2k & \leq \int_{-\pi}^{\pi} \left( 1 + |d w (x,1) |^2 \right)^{1/2} \, dx  \, .
    \label{e:grr5}
\end{align}
\end{Pro}

\begin{figure}[htbp]
    \setlength{\captionindent}{4pt}
    \begin{minipage}[t]{0.5\textwidth}
    \centering\input{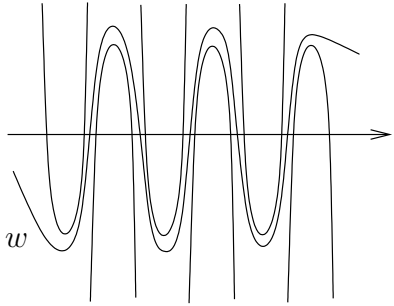}
    \caption{Alternating rescaled grim  reapers.}\label{f:f4}
    \end{minipage}
\end{figure}

\begin{proof}
 For $-k \leq j \leq k$, define translating solutions $u_j$ on
  $j \pi/k < x < (j+1) \pi / k$ by
 \begin{equation}   \label{e:ujxt}
    u_j (x,t)=  (-1)^j \, \left[ u^{k}(x - j\pi / k  ,t) - 2 \, k \right] \,
    ,
 \end{equation}
 where $u^k$ is the scaled grim reaper.
The solutions given by \eqr{e:ujxt}, which alternate between
translating up and down, will be used as barriers; see figure \ref{f:f4}.  
As in the
proof of Proposition \ref{p:grr}, we can choose a compactly
supported function $w(\cdot,0) : \RR \to \RR$ satisfying
\eqr{e:grr4} which is below the upward translating solutions and
above the downward translating solutions.  Combining the existence
results of \cite{EH1} or \cite{EH2} with the maximum principle as
before gives a solution $w(x,t)$ with
\begin{align}
    w( (j+1/2) \pi / k  , 1) < u_j ( (j+1/2) \pi / k  , 1) = - k  &{\text{ for }} j {\text{ even}}, \\
    k = u_j ( (j+1/2) \pi / k  , 1) < w( (j+1/2) \pi / k  , 1)
        &{\text{ for }} j {\text{ odd}}.
\end{align}
The lower bound on length in \eqr{e:grr5} follows immediately.
\end{proof}

We will prove Theorem \ref{t:areabound} by showing that
 the (weighted) area of the graph satisfies
 a differential
inequality which will imply the desired bound (see Lemma
\ref{l:di}).

 We begin with an elementary
 area bound for the graph of a general function $w$:

\begin{Lem}     \label{l:graph}
If $w$, $\phi : \RR^n \to \RR$ are functions  and $\phi$ has
compact support, then
\begin{equation}    \label{e:a2}
   \int \phi^2 (1 + |dw|^2)^{1/2}   \, dx \leq
     \int  \phi^2  \, dx  +  \|w\|_{\infty}  \int |d\phi^2|  \, dx
    + \|w\|_{\infty} \, \int \phi^2  \, |H|  \, dx \, ,
\end{equation}
where $H = - \dv \left( \frac{dw}{(1 + |dw|^2)^{1/2} } \right)$ is
the mean curvature of the graph of $w$.
\end{Lem}

\begin{proof}
 Applying Stokes
theorem to $\dv \left( \frac{ \phi^2 w dw}{(1 + |dw|^2)^{1/2}}
\right)$ gives
\begin{equation}    \label{e:a1}
    \int \phi^2 \frac{ |dw|^2 }{(1 + |dw|^2)^{1/2}} \, dx
    \leq  \int  |d\phi^2| \, \frac{ |w| \, |dw| }
{(1 + |dw|^2)^{1/2}}  \, dx
    + \int \phi^2 |w| \, |H| \, dx  \, .
\end{equation}
Adding $\int \phi^2 \, dx$ to each side gives \eqr{e:a2}.
\end{proof}

When the graph of $w$ is minimal (i.e., $H=0$),
 Lemma \ref{l:graph} gives the well-known area bound
$C \, r^n \, ( 1 + r^{-1} \, \| w \|_{\infty})$.  This linear
dependence on $\| w \|_{\infty}$ is easily seen to be sharp.

\begin{Lem}     \label{l:di}
If $f (t) \geq 0$ and $f^2 \leq - a \, f' + b$ with $a , b
>0$, then $f(T) \leq    \sqrt{2b} + 2a / T$.
\end{Lem}

\begin{proof}
 If $f^2 (T) < 2b$, then we are done.  If
$ f^2 (t) \geq 2b$,
 then $f^2(t)  \leq - 2 a \, f'(t) $ so
\begin{equation}
    \left( 1/f \right)' (t) = - f'(t) / f^2(t)  \geq 1/ (2 a)   \, .
\end{equation}
In particular, if $ f^2 (t) \geq 2b$ on $[t_0, T]$, then
\begin{equation}    \label{e:t0}
    1/ [f(T) ] \geq 1/ [f(t_0)] + (T-t_0)/2a \, .
\end{equation}
We consider two cases.  First, if $ f^2 (t) > 2b$ on $[0, T]$,
then \eqr{e:t0} yields $f(T) \leq 2a / T$.  Otherwise, if $f(t_0)
= \sqrt{2b}$ for some $t_0 < T$, then \eqr{e:t0} gives $f(T) \leq
f(t_0) = \sqrt{2b}$.
\end{proof}

\begin{proof}
(of Theorem \ref{t:areabound}).  By scaling, we can assume that
$r=1$. Within this proof, we write $\| u \|_{\infty}$ for the
$L^{\infty}$ norm of $u$ on $B_1 \times [0,1]$.  Set $\eta (x) =
\max \{(1-|x|), 0 \}$ and define
\begin{equation}
    f(t) = \int \eta^4 \, (1 + |du|^2 )^{1/2} (x,t) \, dx \, .
\end{equation}
(We will omit the $(x,t)$ below.)  Differentiating $f(t)$ and
using Stokes theorem gives
\begin{equation}    \label{e:diffA}
    f'(t) = \int  \eta^4 \, \frac{\langle du , du_t \rangle }{(1 + |du|^2)^{1/2}} \, dx
    = - \int \eta^4 H^2 \, (1 + |du|^2)^{1/2} \, dx
    + 4 \int H \eta^3 \langle d \eta , du \rangle   \, dx \, .
\end{equation}
The absorbing inequality
   $4 |H| \, \eta^3   \leq H^2 \eta^4 / 2
    +  8 \eta^2 $
  then   gives
\begin{equation}    \label{e:diffA2}
   \int \eta^4 H^2 \, (1 + |du|^2)^{1/2} \, dx \leq - 2 f'(t)
    +  16 \int \eta^2 (1 + |du|^2)^{1/2} \, dx \, .
\end{equation}
Applying Lemma \ref{l:graph} with $\phi = \eta$ and using an
absorbing inequality gives
\begin{equation}    \label{e:diffA3}
    16 \, \int \eta^2 (1 + |du|^2)^{1/2} \, dx
\leq C_1 \, (1 + \|u\|_{\infty}) +  C_2 \,
    \|u\|_{\infty}^2 + 1/2
    \, \int \eta^4 \, H^2 \, dx
    \, .
\end{equation}
Combining \eqr{e:diffA2} and \eqr{e:diffA3} gives
\begin{equation}    \label{e:diffA4}
   \int \eta^4 H^2 \, (1 + |du|^2)^{1/2} \, dx \leq - 4 f'(t)
    +  C_3 ( 1 + \|u\|_{\infty}^2)  \, .
\end{equation}
 After applying Lemma \ref{l:graph} with $\phi = \eta^2$, the
Cauchy-Schwarz inequality and  \eqr{e:diffA4} give
\begin{align}        \label{e:A}
    f^2(t) &
    \leq  C_4  \left(1 + \|u\|_{\infty}^2  +    \|u\|_{\infty}^2
    \, \left( \int \eta^4 \, H^2 \, dx \right) \, \left( \int \eta^4 \,
    dx\right)
    \right)
    \notag  \\
    &
    \leq  C_4 (1 + \|u\|_{\infty}^2)
+ C_5 \, \|u\|_{\infty}^2 \left( - 4 f'(t)
    +  C_3 ( 1 + \|u\|_{\infty}^2) \right)  \, .
\end{align}
Finally, applying Lemma \ref{l:di} gives the theorem since, by  Lemma
\ref{l:sphere},  $\|u\|_{\infty} \leq \sqrt{2n+1} + \sup_{
B_{\sqrt{2n+1} } } |u (\cdot , 0) |$.
\end{proof}


\begin{thebibliography}{999}


\bibitem[BDM]{BDM}
E. Bombieri, E. De Giorgi, and M. Miranda,  Una maggiorazione a
priori relativa alla ipersuperfici minimali non parametriche, {\it
Arch. Rational Mech. Anal.}, 32 (1969) 255--267.
\bibitem[E]{E}
K. Ecker,  Lectures on regularity for mean curvature flow,
preprint.
\bibitem[EH1]{EH1}
K. Ecker and G. Huisken,  Mean curvature evolution of entire
graphs, {\it Annals of Math.}, 130 (1989) 453--471.
\bibitem[EH2]{EH2}
K. Ecker and G. Huisken,  Interior estimates for hypersurfaces
moving by mean curvature, {\it Invent. Math.}, 105 (1991)
547--569.
\bibitem[F1]{F1}
R. Finn,  On equations of minimal surface type, {\it Annals of
Math.}, 60 (1954) 397--416.
\bibitem[F2]{F2}
R. Finn,  Remarks relevant to minimal surfaces, and to surfaces of
prescribed mean curvature, {\it J. Analyse Math.}, 14 (1965)
265--296.
\bibitem[GiTr]{GiTr}
D. Gilbarg and N. Trudinger, Elliptic partial differential
equations of second order, Springer Verlag  (1983).
\bibitem[K]{K}
N. Korevaar, An easy proof of the interior gradient bound for
solutions to the prescribed mean curvature equation, Nonlinear
functional analysis and its applications, Part 2, {\it Proc. of
Symposia in Pure  Math.} 45 (1986)  81--89.


\end{thebibliography}
\end{document}